# PRISMATIC LOGARITHM AND PRISMATIC HOCHSCHILD HOMOLOGY VIA NORM

ZHOUHANG MAO

ABSTRACT. In this brief note, we present an elementary construction of the first Chern class of Hodge–Tate crystals in line bundles using a refinement of the prismatic logarithm, which should be comparable to the one considered by Bhargav BHATT. The key to constructing this refinement is Yuri SULYMA's norm on (animated) prisms. We explain the relation of this construction to prismatic Witt vectors, as a generalization of Kaledin's polynomial Witt vectors. We also propose the prismatic Hochschild homology as a noncommutative analogue of prismatic de Rham complex.

## 1. INTRODUCTION

Let $X$ a bounded $p$-adic formal scheme, and $\mathcal{E}$ a vector bundle on $X$. In [BL22, §7 & §9.2], Bhatt–Lurie constructed the prismatic[1.1] Chern classes $c_i^{\mathbb{A}}(\mathcal{L}) \in H_{\mathbb{A}}^{2i}(X)\{i\}$, using the prismatic logarithm

$$\log_{\mathbb{A}}^{\mathrm{BL}} : T_p(\mathbb{G}_m(\overline{A})) \longrightarrow A\{1\}$$

on the Tate module for bounded prisms $(A, I)$ with $\overline{A} := A/I$. In [Bha23, Cons 7.6], Bhatt observed that the prismatic Chern classes $c_i^{\mathbb{A}}(\mathcal{E})$ only depend on the pullback of the vector bundle $\mathcal{E}$ to the Hodge–Tate stack $X^{\mathrm{HT}}$, and he defined the prismatic Chern classes $c_i(E) \in H_{\mathbb{A}}^{2i}(X)\{i\}$ for Hodge–Tate crystals $E$ in vector bundles, obtained by pulling back the generators of the prismatic cohomology of the stack $\mathrm{BGL}_n$.

In this brief note, we propose an elementary construction of the first Chern class $c_1(L)$ of Hodge–Tate crystals $L$ in line bundles when $p > 2$ is an odd prime, using a variant

$$\mathrm{d}\log_{\mathbb{A}} : \mathbb{G}_m(\overline{A}) \longrightarrow A\{1\}[1]$$

of the prismatic logarithm, for any animated prism $(A, I)$ with $\overline{A} := A/I$. The key to this construction is to produce Sulyma's norm maps [Sul23, §3.1] on animated prisms. We show that this refines Bhatt–Lurie's prismatic logarithm $\log_{\mathbb{A}}^{\mathrm{BL}}$.

This construction comes out of discussions with Alexander PETROV who expects some relation between Kaledin's polynomial Witt vectors and crystalline Chern classes. Let $(A, I)$ be a transversal prism. Using relative HHR norms, we propose a prismatic analogue of Kaledin's polynomial Witt vectors, called *prismatic Witt vectors*. More precisely, for every positive integer $r \in \mathbb{N}_{>0}$, the $r$-truncated prismatic Witt vectors functor

$$\mathbb{A}_r : D(\overline{A}) \longrightarrow D(A/I_r)$$

---

[1.1]. In fact, they refine it to the syntomic Chern classes. We will not consider this refinement in this brief note.

is a polynomial functor of degree $p^{r-1}$, which is compatible with base change. It is an analogue since Kaledin's polynomial Witt vectors, as well as their generalizations in [DKNP23, Rea23], can be constructed out of norms, which will be explained in our forthcoming work [Mao]. As a consequence, these functors give rise to a family of polynomial functors sending Hodge–Tate crystals to prismatic crystals. This construction is closely related to our refined logarithm. Moreover, when $(A, I)$ is a transversal perfect prism, the functor $\mathbb{A}_r$ coincides with the $p$-typical polynomial Witt vector functor $W_{r,p}(\overline{A}; -)$ on finite free $\overline{A}$-modules (viewed as symmetric $\overline{A}$-bimodules).

In [BMS19], the authors related topological Hochschild homology and its variants to the absolute prismatic cohomology. It is a natural question whether we can relate some modified Hochschild homology to relative prismatic cohomology, which suggests the former being noncommutative relative prismatic cohomology. [BMS19, §11], the authors gave a positive answer to this when the base prism $(A, I)$ is the Breuil–Kisin prism. In Section 3, based on prismatic Witt vectors, we propose a candidate — the *prismatic Hochschild homology* $\text{HH}^{\mathbb{A}}(R/A)$ for every associative $A/I$-algebra $R$ (even for every $A/I$-linear DG-category). We conjecture an HKR type theorem for $p$-completely smooth $A/I$-algebra $R$, relating the $r$-truncated prismatic Hochschild homology groups $\pi_n\big(\text{HH}^{\mathbb{A}}(R/A)^{C_{p^{r-1}}}\big)$ to the relative prismatic cohomology $\mathbb{A}_{R/A}$. When $r = 1$, this reduces to the usual HKR theorem. We also sketch a proof of it for $p$-completed polynomial $A/I$-algebras.

*Acknowledgments.* We thank Yuri Sulyma for explaining and discussions of his Tambara functors associated to prisms, Wolfgang Steimle and Kaif Hilman for letting us know the concept of the Hill–Hopkins–Ravenel norm and its relative version, and Alexander Petrov for explaining the potential relation between Kaledin's polynomial Witt vectors and crystalline Chern classes, and a proof in the special case of the first crystalline Chern class. We would also like to thank Bastiaan Cnossen, Matthew Morrow, and Maxime Ramzi.

## 2. Refined logarithm

In this section, we first review Yuri Sulyma's norm maps for transveral prisms $(A, I)$. Then we briefly explain Bhatt–Lurie's prismatic logarithm in [BL22, §2], construct our refined prismatic logarithm, and compare it with Bhatt–Lurie's on transversal prisms. Finally, we briefly indicate how to extend the constructions to animated prisms.

2.1. **Norms on transversal prisms.** We recollect Yuri Sulyma's norm on transversal prisms. Let $(A, I)$ be a transversal prism (that is to say, the quotient ring $A/I$ is $p$-torsion-free). We will denote by $\overline{A}$ the quotient ring $A/I$, by $I_r$ the invertible ideal $I\varphi^*(I)\cdots(\varphi^{r-1})^*(I) \subseteq A$ (when $r=0$, it is the unit ideal) (cf. [BL22, Nota 2.2.2]), and by $J_r$ the invertible ideal $I^p\big(\varphi^*(I)\cdots(\varphi^{r-1})^*(I)\big)^{p-1} = I \cdot I_r^{p-1}$ (when $r=0$, it is $I$). We review [BL22, Lem 2.2.8] by strengthening it a bit:

**Lemma 2.1.** *Let $(A, I)$ be a transversal prism. For every $r \in \mathbb{N}_{>0}$, the canonical map*

$$f: (\varphi^r)^*(I)/J_r(\varphi^r)^*(I) \longrightarrow A/J_r$$

*is a monomorphism, whose image is the principal ideal $(p) \subseteq A/J_r$.*

*Proof.* We adapt the proof of [BL22, Lem 2.2.8]. First, in the quotient ring $A/p$, the image of $(\varphi^r)^*(I)$ coincides with the image of $I^{p^r}$, which subsequently coincides with the image of $J_r$. This implies[2.1] that $\mathrm{Im}(f) \subseteq (p)$.

The proof of [BL22, Lem 2.2.5] implies that the quotient ring $A/J_r$ is $p$-torsion-free, thus there exists a unique map $f_0 \colon (\varphi^r)^*(I)/J_r(\varphi^r)^*(I) \to A/J_r$ of $A$-modules such that $f = p f_0$. We now show that the map $f_0$ is an isomorphism, or equivalently, it is surjective. For every $x \in I$, we have

$$\begin{aligned}
\varphi^r(x) &= \varphi^{r-1}(x^p + p\,\delta(x)) \\
&= \varphi^{r-1}(x)^{p-1}\,\varphi^{r-1}(x) + p\,\varphi^{r-1}(\delta(x)) \\
&\in \varphi^{r-1}(x)^{p-1}\,(J_{r-1} + (p)) + p\,\varphi^{r-1}(\delta(x)) \\
&\subseteq J_r + p\,(\varphi^{r-1})^*(I)^{p-1} + p\,\varphi^{r-1}(\delta(x)),
\end{aligned}$$

where the first containment uses the fact that $\mathrm{Im}(f) \subseteq (p)$ for $r-1$. It follows that $f_0(x) \equiv \varphi^{r-1}(\delta(x)) \pmod{(\varphi^{r-1})^*(I)}$. Since $(A,I)$ is a prism, the elements $\{\delta(x) \,|\, x \in I\}$ generate the unit ideal of $A/(\varphi^{r-1})^*(I)$, and hence also in $A/J_r$. □

**Corollary 2.2. (cf. [Sul23, Lem 3.1])** *Let $(A, I)$ be a transversal prism. Then the commutative diagram*

$$\begin{array}{ccc}
A/J_r(\varphi^r)^*(I) & \longrightarrow & A/(\varphi^r)^*(I) \\
\downarrow & \square & \downarrow \\
A/J_r & \longrightarrow & A/(J_r, (\varphi^r)^*(I))
\end{array}$$

*of commutative $A$-algebras is Cartesian in $D(A)$. The same holds if we replace $J_r$ by $I_r^j$ for $0 \leq j \leq p-1$.*

Now we construct the norm as follows.

*Construction* 2.3. **([Sul23, Cons 3.3])** Let $(A, I)$ be a transversal prism, and $r \in \mathbb{N}_{>0}$. Note that there is a multiplicative map

$$\begin{aligned}
A/I_r &\longrightarrow A/I_r \times_{A/(I_r, (\varphi^r)^*(I))} A/(\varphi^r)^*(I) \\
x &\longmapsto (x^p, \varphi(x))
\end{aligned}$$

which is polynomial of degree $\leq p$ (we tacitly invoked Lemma 2.1). Thanks to Corollary 2.2 (for $I_r$), we obtain a multiplicative map

$$N_r \colon A/I_r \longrightarrow A/I_{r+1}$$

which is polynomial of degree $\leq p$ as well. In particular, this gives rise to a map

$$N_r \colon \mathbb{G}_m(A/I_r) \longrightarrow \mathbb{G}_m(A/I_{r+1})$$

of abelian groups.

*Remark* 2.4. Let $(A, I)$ be a transversal prism. In [Sul23, Cons 3.3], the author constructed a $\mathbb{T}$-Tambara functor $\underline{(A, I)}$, such that

$$\underline{(A, I)}^{C_{p^{r-1}}} = A/I_r$$

---

2.1. Note that, the proof of this inclusion only uses the fact that $(A, I)$ is a pre-prism, which is generalized in Construction 2.28.

(a fortiori with trivial Weyl action) for $r \in \mathbb{N}_{>0}$. The norm maps $\underline{(A,I)}^{C_{p^{r-1}}} \to \underline{(A,I)}^{C_{p^r}}$ are given by the maps $N_r$ in Construction 2.3, the restriction maps $\underline{(A,I)}^{C_{p^r}} \to \underline{(A,I)}^{C_{p^{r-1}}}$ are given by the canonical quotient maps $A/I_{r+1} \to A/I_r$, and the transfer maps $\underline{(A,I)}^{C_{p^{r-1}}} \to \underline{(A,I)}^{C_{p^r}}$ are given by

$$A/I_r \longrightarrow A/I_r \times_{A/(I_r, (\varphi^r)^*(I))} A/(\varphi^r)^*(I)$$
$$x \longmapsto (p\,x, 0),$$

invoking Corollary 2.2.

*Example* 2.5. ([Sul23, Thm 3.19]) Let $(A,I)$ be a transversal perfect prism. Then the $\mathbb{T}$-Tambara functor $\underline{(A,I)}$ coincides with the $\mathbb{T}$-Tambara functor $\underline{W}(\overline{A})$ of Witt vectors.

*Question* 1. How to extend the construction in [Sul23, Cons 3.3] (reviewed in Remark 2.4) to derived prisms, giving rise to a functor from derived prisms to derived $\mathbb{T}$-Tambara functors?

2.2. **Refining $r$-truncated logarithm.** Bhatt–Lurie's prismatic logarithm $\log_{\mathbb{A}}^{\mathrm{BL}} \colon T_p(\mathbb{G}_m(\overline{A})) \to A\{1\}$ is constructed out of its "$r$-truncated" versions, being the composite

$$\log_{\mathbb{A}}^{(r),\mathrm{BL}} \colon T_p(\mathbb{G}_m(\overline{A})) \longrightarrow (1+I)_{\mathrm{rk}=1} \xrightarrow{(-)^{p^{r-1}}-1} I_r \longrightarrow I_r/I_r^2.$$

In this subsection, we refine these maps.

*Construction* 2.6. Let $R$ be a commutative ring. Then we have a map

$$\mathrm{d}\log \colon \mathbb{G}_m(R) \longrightarrow \Omega^1_{R/\mathbb{Z}}$$
$$f \longmapsto f^{-1}\mathrm{d}f$$

of abelian groups. Left deriving this map, we get a map $\mathrm{d}\log \colon \mathbb{G}_m(R) \to L_{R/\mathbb{Z}}$ for any animated ring $R$.

*Construction* 2.7. Let $A$ be a commutative ring, and $I \subseteq A$ a square-zero ideal. Then we have a short exact sequence

$$0 \longrightarrow I \xrightarrow{\exp} \mathbb{G}_m(A) \longrightarrow \mathbb{G}_m(A/I) \longrightarrow 1$$

where the first map $\exp \colon I \to \mathbb{G}_m(A)$ is given by $x \mapsto 1+x$ for $x \in I \subseteq A$. Consequently, we get a boundary map $\mathbb{G}_m(A/I) \to I[1]$ in $D_{\geqslant 0}(\mathbb{Z})$, which is understood as a shifted version of logarithm (or $y \mapsto y-1$). Left derive this construction from $(A,I) = (\mathbb{Z}[x_i, y_j]/(y_k y_l), (y_j))$, we get a map $\mathbb{G}_m(A/I) \to I[1]$ in $D_{\geqslant 0}(\mathbb{Z})$ for any animated square-zero ring-ideal pairs $(A,I)$.

In particular, given an animated ring-ideal pair $(A,I)$, we can form an animated square-zero ring-ideal pair $(A/I^2, I/I^2)$, and in this case, the map above is denoted by $\mathrm{d}\log_{(A,I)} \colon \mathbb{G}_m(A/I) \to (I/I^2)[1]$.

*Remark* 2.8. The notation $\mathrm{d}\log$ in Construction 2.7 is justified by the fact that $\mathrm{d}\log_{(A,I)}$ is canonically identified with the composite

$$\mathbb{G}_m(A/I) \xrightarrow{\mathrm{d}\log} L_{(A/I)/\mathbb{Z}} \longrightarrow L_{(A/I)/A} \simeq (I/I^2)[1]$$

for any animated ring-ideal pair $(A,I)$. Indeed, it suffices to produce this map *functorially* on the set $\{(\mathbb{Z}[x_i, y_j], (x_i))\}$ of compact projective generators (cf. [Mao21, Thm 1.1]). In this case, the transitivity sequence

$$L_{(A/I)/A}[-1] \longrightarrow L_{A/\mathbb{Z}} \otimes_A^{\mathbb{L}} (A/I) \longrightarrow L_{(A/I)/\mathbb{Z}}$$

can be identified with the conormal sequence

$$0 \longrightarrow I/I^2 \xrightarrow{d} \Omega^1_{A/\mathbb{Z}} \otimes_A (A/I) \longrightarrow \Omega^1_{(A/I)/\mathbb{Z}} \longrightarrow 0$$

which fits into a morphism

$$\begin{array}{ccccccccc}
0 & \longrightarrow & I/I^2 & \xrightarrow{\exp} & \mathbb{G}_m(A/I^2) & \longrightarrow & \mathbb{G}_m(A/I) & \longrightarrow & 0 \\
& & \| & & \downarrow & & \downarrow & & \\
0 & \longrightarrow & I/I^2 & \xrightarrow{d} & \Omega^1_{A/\mathbb{Z}} \otimes_A (A/I) & \longrightarrow & \Omega^1_{(A/I)/\mathbb{Z}} & \longrightarrow & 0
\end{array},$$

of short exact sequences, where the middle vertical map is $d\log \colon \mathbb{G}_m(A/I^2) \to \Omega^1_{A/\mathbb{Z}} \otimes_A (A/I), [f] \mapsto df \otimes [f]^{-1}$, and the right vertical map is $d\log \colon \mathbb{G}_m(A/I) \to \Omega^1_{(A/I)/\mathbb{Z}}$ in Construction 2.6. This realizes the extension class $\mathbb{G}_m(A/I) \to (I/I^2)[1]$ as the composite $\mathbb{G}_m(A/I) \xrightarrow{d\log} L_{(A/I)/\mathbb{Z}} = \Omega^1_{(A/I)/\mathbb{Z}} \to (I/I^2)[1] \simeq L_{(A/I)/A}$, functorially in ring ideal pairs $(A, I)$ with both $A$ and $A/I$ being polynomial rings. Left deriving this, we get the result for animated ring-ideal pairs $(A, I)$.

*Remark* 2.9. (M. RAMZI) We can replace $\mathbb{G}_m$ by $\mathrm{GL}_1$ in Construction 2.7, and the map for $\mathbb{G}_m$ can be obtained by composing with the canonical map $\mathbb{G}_m \to \mathrm{GL}_1$. In fact, this construction works $\mathbb{E}_\infty$-ly. More precisely, let $(A, I)$ be a connective square-zero $\mathbb{E}_\infty$-ring-ideal pair. Then we can construct a fiber sequence

$$I[1] \longrightarrow \mathrm{pic}(A) \longrightarrow \mathrm{pic}(A/I)$$

of spectra, where $\mathrm{pic}(-)$ is the Picard spectrum, as follows. First, we realize the square-zero extension $A$ of $A/I$ by $I$ as a pullback

$$\begin{array}{ccc}
A & \longrightarrow & A/I \\
\downarrow & & \downarrow \\
A/I & \longrightarrow & A/I \oplus I[1]
\end{array}$$

by definition, where the bottom map is associated to the zero derivation, while the right vertical map is part of data of the square-zero extension. This gives rise to a Cartesian square

$$\begin{array}{ccc}
\mathrm{pic}(A) & \longrightarrow & \mathrm{pic}(A/I) \\
\downarrow & & \downarrow \\
\mathrm{pic}(A/I) & \longrightarrow & \mathrm{pic}(A/I \oplus I[1])
\end{array}$$

of spectra (cf. [Lur18, Prop 16.2.2.1]). The map $A/I \to A/I \oplus I[1]$ of $\mathbb{E}_\infty$-rings fits into a Cartesian square

$$\begin{array}{ccc}
\mathbb{S} & \longrightarrow & \mathbb{S} \oplus I[1] \\
\downarrow & & \downarrow \\
A/I & \longrightarrow & A/I \oplus I[1]
\end{array}$$

of connective $\mathbb{E}_\infty$-rings. Taking $\mathrm{pic}(-)$ (here we need the connectivity of $I$, not just $(-1)$-connectivity), we get a Cartesian square

$$\begin{array}{ccc}
\mathrm{pic}(\mathbb{S}) & \longrightarrow & \mathrm{pic}(\mathbb{S} \oplus I[1]) \\
\downarrow & & \downarrow \\
\mathrm{pic}(A/I) & \longrightarrow & \mathrm{pic}(A/I \oplus I[1])
\end{array}$$

of spectra. Now we have

$$\begin{aligned}\mathrm{fib}(\mathrm{pic}(A)\to\mathrm{pic}(A/I)) &\simeq \mathrm{fib}(\mathrm{pic}(A/I)\to\mathrm{pic}(A/I\oplus I[1]))\\ &\simeq \mathrm{fib}(\mathrm{pic}(\mathbb{S})\to\mathrm{pic}(\mathbb{S}\oplus I[1]))\\ &\simeq \Omega\,\mathrm{fib}(\mathrm{pic}(\mathbb{S}\oplus I[1])\to\mathrm{pic}(\mathbb{S}))\\ &\simeq \mathrm{fib}(\mathrm{GL}_1(\mathbb{S}\oplus I[1])\to\mathrm{GL}_1(\mathbb{S}))\\ &\simeq I[1],\end{aligned}$$

where we use the fact that $\Omega\circ\mathrm{pic}=\mathrm{GL}_1$.

*Construction* 2.10. Let $(A,I)$ be a transversal prism. Then we define the map $\mathrm{d}\log_{\mathbb{A}}^{(r)}$ to be the composite map

$$\mathbb{G}_m(\overline{A})\xrightarrow{N_{r-1}\circ\cdots\circ N_1}\mathbb{G}_m(A/I_r)\xrightarrow{\mathrm{d}\log_{(A,I_r)}}(I_r/I_r^2)[1]$$

in $D(\mathbb{Z})$, where the first map is in Construction 2.3, and the second map is in Construction 2.7. The map $\mathrm{d}\log_{\mathbb{A}}^{(r)}$ is informally understood as a deformation of $y\mapsto y^{p^{r-1}}-1$. Since the object $(I_r/I_r^2)[1]\in D(\mathbb{Z})$ is $p$-complete, we can identify the $p$-completion of the map $\mathrm{d}\log_{\mathbb{A}}^{(r)}$ as a map

$$\widehat{\mathrm{d}\log}_{\mathbb{A}}^{(r)}:\mathbb{G}_m(\overline{A})_p^{\wedge}\longrightarrow (I_r/I_r^2)[1]$$

in $D(\mathbb{Z})$.

*Remark* 2.11. Let $(A,I)$ be a transveral perfect prism. Then the iterated norm map

$$N_{r-1}\circ\cdots\circ N_1:\overline{A}\longrightarrow A/I_r$$

coincides with the Teichmüller map $\overline{A}\to W_r(\overline{A})$ under the identification Example 2.5, thus the iterated norm map is understood as a *prismatic Teichmüller map*.

We now show that the map $\widehat{\mathrm{d}\log}_{\mathbb{A}}^{(r)}:\mathbb{G}_m(\overline{A})_p^{\wedge}\to(I_r/I_r^2)[1]$ refines Bhatt–Lurie's map $\log_{\mathbb{A}}^{(r),\mathrm{BL}}:T_p(\mathbb{G}_m(\overline{A}))\to(I_r/I_r^2)[1]$. The key is Remark 2.14, which tells us how to identify the first homology group $H_1(-)$ of refined logarithm in a more general context.

*Remark* 2.12. Let $M$ be an abelian group. Then its derived $p$-completion $M_p^{\wedge}$ can be identified with the cofiber of the composite map

$$\prod_{i\in\mathbb{N}_{>0}}M\xrightarrow{(0,-)}\prod_{i\in\mathbb{N}}M\xrightarrow{p\mathrm{shift}-\mathrm{id}}\prod_{i\in\mathbb{N}}M, \qquad(2.1)$$

and the map $M_p^{\wedge}\to M/^{\mathbb{L}}p^n$ can be identfied with the vertical map between horizontal cofibers of the commutative diagram

$$\begin{array}{ccc}\prod_{i\in\mathbb{N}_{>0}}M & \longrightarrow & \prod_{i\in\mathbb{N}}M\\ \downarrow & & \downarrow\\ M & \xrightarrow{p^n} & M\end{array}$$

where the left vertical map is the $n$-th projection $(x_i)_{i\in\mathbb{N}_{>0}}\mapsto x_n$, and the right vertical map is $(y_i)\mapsto y_1+p\,y_2+\cdots+p^n y_n$.

**Notation 2.13.** *Let $R$ be a commutative ring. We will denote by $R^{\flat}$ the sequential limit*

$$\lim\left(\cdots\xrightarrow{(-)^p}R\xrightarrow{(-)^p}R\right)$$

of multiplicative monoids, and by $R^{\flat\times}$ the invertible elements of $R^\flat$.

*Remark* 2.14. Let $A$ be a commutative ring, and $I \subseteq A$ an ideal. Then the first homology group $H_1(-)$ of the $p$-completed map

$$(\mathrm{d}\log_{(A,I)})_p^\wedge : \mathbb{G}_m(A/I)_p^\wedge \longrightarrow (I/I^2)_p^\wedge[1]$$

can be identified with the boundary map of the snake long exact sequence associated to the morphism

$$\begin{array}{ccccccccc}
0 & \longrightarrow & \prod_{i\in\mathbb{N}_{>0}} I/I^2 & \longrightarrow & \prod_{i\in\mathbb{N}_{>0}} \mathbb{G}_m(A/I^2) & \longrightarrow & \prod_{i\in\mathbb{N}_{>0}} \mathbb{G}_m(A/I) & \longrightarrow & 1 \\
& & \downarrow & & \downarrow & & \downarrow & & \\
0 & \longrightarrow & \prod_{i\in\mathbb{N}} I/I^2 & \longrightarrow & \prod_{i\in\mathbb{N}} \mathbb{G}_m(A/I^2) & \longrightarrow & \prod_{i\in\mathbb{N}} \mathbb{G}_m(A/I) & \longrightarrow & 1
\end{array}$$

of short exact sequences, where vertical maps are of the form (2.1). Concretely, given an element $(x_i)_{i\in\mathbb{N}} \in T_p(\mathbb{G}_m(A/I))$ with $x_0 = 1$, we pick an $A/I^2$-lift $(y_i)_{i\in\mathbb{N}} \in \prod_{i\in\mathbb{N}} \mathbb{G}_m(A/I^2)$ with $y_0 = 1$, and the image of $(x_i)_{i\in\mathbb{N}} \in T_p(\mathbb{G}_m(A/I))$ under the boundary map is given by

$$\left(y_i^{p^i} - 1\right)_{i\in\mathbb{N}} \in \varprojlim_{i\in\mathbb{N}} (I/I^2)/p^i \subseteq \prod_{i\in\mathbb{N}} (I/I^2)/p^i.$$

In particular, if there is an $A/I^2$-lift $(y_i)_{i\in\mathbb{N}} \in (A/I^2)^{\flat\times}$, then the image of $(x_i) \in T_p(\mathbb{G}_m(A/I))$ is simply given by $y_1 - 1 \in I/I^2$.

In our applications, the canonical lift $(y_i)_{i\in\mathbb{N}} \in (A/I^2)^{\flat\times}$ for $(x_i)_{i\in\mathbb{N}} \in T_p(\mathbb{G}_m(A/I))$ does exist:

**Lemma 2.15.** *Let $(A, I)$ be a prism, and $J \subseteq \mathrm{Rad}(A)$ an ideal lying in the Jacobson radical which contains $I$. Then the quotient map $A \twoheadrightarrow A/J$ induces an isomorphism $A^{\flat\times} \to (A/J)^{\flat\times}$ of abelian groups.*

*Proof.* Since the ideal $J \subseteq A$ lies in the Jacobson radical of $A$, it follows that the map in question is surjective. Now we consider the composable maps $A \twoheadrightarrow A/J \twoheadrightarrow A/I$ of rings, which induces composable surjective maps $A^{\flat\times} \to (A/J)^{\flat\times} \to (A/I)^{\flat\times}$, whose composite is an isomorphism by [BL22, Prop 2.7.3]. The result then follows. $\square$

**Proposition 2.16.** *Let $(A, I)$ be a transversal prism, and $r \in \mathbb{N}_{>0}$. Then the first homology group $H_1(-)$ of the map*

$$\widehat{\mathrm{d}\log}_{\mathbb{A}}^{(r)} : \mathbb{G}_m(\overline{A})_p^\wedge \longrightarrow (I_r/I_r^2)[1]$$

*in $D(\mathbb{Z})$ can be identified with Bhatt–Lurie's map $\log_{\mathbb{A}}^{(r),\mathrm{BL}}$.*

*Proof.* Applying Remark 2.14 and Lemma 2.15 to $(A, I_r)$. $\square$

*Remark* 2.17. We are not aware of the best comparison between our refined $r$-truncated prismatic logarithm with the $r$-truncated prismatic logarithm $(1+I)_{\mathrm{rk}=1} \to I_r/I_r^2$. However, the prismatic logarithm $(1+I)_{\mathrm{rk}=1} \to I_r/I_r^2$ can be related to $\log_{(A,I_r)}$ as follows. We examine the morphism

$$\begin{array}{ccccccccc}
0 & \longrightarrow & I_r/I_r^2 & \longrightarrow & \mathbb{G}_m(A/I_r^2) & \longrightarrow & \mathbb{G}_m(A/I_r) & \longrightarrow & 1 \\
& & \downarrow & & \downarrow & & \downarrow & & \\
0 & \longrightarrow & I_r/I_r^2 & \longrightarrow & \mathbb{G}_m(A/I_r^2) & \longrightarrow & \mathbb{G}_m(A/I_r) & \longrightarrow & 1
\end{array}$$

of short exact sequences, where vertical maps are given by $p^{r-1}$-scaling. The boundary map gives rise to a map $\mu_{p^{r-1}}(A/I_r) \to (I_r/I_r^2)/p^{r-1}$, and the composite map $(1+I)_{\mathrm{rk}=1} \xrightarrow{\mathrm{id}} \mu_{p^{r-1}}(A/I_r) \to (I_r/I_r^2)/p^{r-1}$ coincides with the prismatic logarithm after modulo $p^{r-1}$. A shift of this can be rewritten as the composite

$$(1+I)_{\mathrm{rk}=1}[1] \longrightarrow \mathbb{G}_m(A/I_r)/^{\mathbb{L}}p^{r-1} \xrightarrow{\log_{(A,I_r)}/^{\mathbb{L}}p^{r-1}} (I_r/I_r^2)[1]/^{\mathbb{L}}p^{r-1}.$$

2.3. **Refining logarithm.** Let $(A, I)$ be a transversal prism. So far, we have constructed $r$-truncated refined logarithm

$$\widetilde{\mathrm{dlog}}^{(r)}_{\mathbb{\Delta}} \colon \mathbb{G}_m(\overline{A}) \longrightarrow (I_r/I_r^2)[1]$$

in $D(\mathbb{Z})$ for every $r \in \mathbb{N}_{>0}$. We now put them together to get a (non-truncated) refined logarithm.

Recall that the *Breuil–Kisin twist* $A\{1\}$ is the sequential limit

$$\cdots \xrightarrow{p^{-1}} I_r/I_r^2 \xrightarrow{p^{-1}} \cdots \xrightarrow{p^{-1}} I_2/I_2^2 \xrightarrow{p^{-1}} I_1/I_1^2 = I/I^2$$

with surjective transition maps being $p^{-1}$ by Lemma 2.1, and Bhatt–Lurie's prismatic logarithm $\log^{\mathrm{BL}}_{\mathbb{\Delta}} \colon T_p(\mathbb{G}_m(\overline{A})) \to A\{1\}$ is obtained by taking the sequential limit along $r \in \mathbb{N}$. In order to pass to the sequential limit along $r \in \mathbb{N}_{>0}$ for our refined logarithm, we have to specify the compatibility data. It suffices to specify a homotopy, for every $r \in \mathbb{N}_{>0}$, making the diagram

$$\begin{array}{ccc} \mathbb{G}_m(A/I_r) & \xrightarrow{\widetilde{\mathrm{dlog}}^{(r)}_{\mathbb{\Delta}}} & (I_r/I_r^2)[1] \\ \downarrow{\scriptstyle N_r} & & \uparrow{\scriptstyle p^{-1}} \\ \mathbb{G}_m(A/I_{r+1}) & \xrightarrow{\widetilde{\mathrm{dlog}}^{(r+1)}_{\mathbb{\Delta}}} & (I_{r+1}/I_{r+1}^2)[1] \end{array}$$

commute. Unfortunately, we could only cook up these homotopies for odd primes $p > 2$. The key is Lemma 2.1, which implies that the image of $(\varphi^r)^*(I)$ in the quotient ring $A/I_r^2$ since $II_r^{p-1} \subseteq I_r^{p-1} \subseteq I_r^2$ as $p - 1 \geq 2$. We first recall that, the map $p^{-1} \colon I_{r+1}/I_{r+1}^2 \to I_r/I_r^2$ factors as

$$I_{r+1}/I_{r+1}^2 \twoheadrightarrow I_{r+1}/I_r I_{r+1} \xrightarrow[\simeq]{p^{-1}} I_r/I_r^2,$$

thus it suffices to produce a homotopy making the diagram

$$\begin{array}{ccc} \mathbb{G}_m(A/I_r) & \xrightarrow{\widetilde{\mathrm{dlog}}^{(r)}_{\mathbb{\Delta}}} & (I_r/I_r^2)[1] \\ \downarrow{\scriptstyle N_r} & & \simeq\uparrow{\scriptstyle p^{-1}} \\ \mathbb{G}_m(A/I_{r+1}) & \longrightarrow & (I_{r+1}/I_r I_{r+1})[1] \end{array}$$

where the bottom horizontal arrow is obtained by applying Construction 2.7 to the square-zero extension $(A/I_r I_{r+1}, I_{r+1}/I_r I_{r+1})$. Since the right vertical map is an equivalence, we could instead specify a homotopy making the diagram

$$\begin{array}{ccc} \mathbb{G}_m(A/I_r) & \xrightarrow{\widetilde{\mathrm{dlog}}^{(r)}_{\mathbb{\Delta}}} & (I_r/I_r^2)[1] \\ \downarrow{\scriptstyle N_r} & & \simeq\downarrow{\scriptstyle p} \\ \mathbb{G}_m(A/I_{r+1}) & \longrightarrow & (I_{r+1}/I_r I_{r+1})[1] \end{array} \qquad (2.2)$$

commute.

*Construction* 2.18. Let $(A, I)$ be a transversal prism, and $p > 2$ an odd prime. Then we construct a map
$$\mathbb{G}_m(A/I_r^2) \longrightarrow \mathbb{G}_m(A/I_r I_{r+1})$$
which fits into a morphism

$$\begin{array}{ccccccccc}
0 & \longrightarrow & I_r/I_r^2 & \longrightarrow & \mathbb{G}_m(A/I_r^2) & \longrightarrow & \mathbb{G}_m(A/I_r) & \longrightarrow & 0 \\
& & \simeq \downarrow p & & \downarrow & & \downarrow N_r & & \\
0 & \longrightarrow & I_{r+1}/I_r I_{r+1} & \longrightarrow & \mathbb{G}_m(A/I_r I_{r+1}) & \longrightarrow & \mathbb{G}_m(A/I_{r+1}) & \longrightarrow & 0
\end{array} \quad (2.3)$$

of short exact sequences. Note that $I_r I_{r+1} = I_r^2 (\varphi^r)^*(I)$, and by Corollary 2.2, it suffices to produce two multiplicative maps $A/I_r^2 \to A/I_r^2$ and $A/I_r^2 \to A/(\varphi^r)^*(I)$ which coincide after passing to the quotient ring $A/(I_r^2, (\varphi^r)^*(I))$. Indeed, the first map is the $p$-th power map $(-)^p$, and the second map is the composite map
$$A/I_r^2 \longrightarrow A/(\varphi^{r-1})^*(I) \xrightarrow{\varphi} A/(\varphi^r)^*(I).$$

The two map coincides after passing to the quotient ring $A/(I_r^2, (\varphi^r)^*(I))$ by Lemma 2.1 and the definition of $\delta$-rings. It is direct to check that this map fits into (2.3).

The extension class of (2.3) in Construction 2.18 gives rise to a homotopy for (2.2), which concludes our construction of refined logarithm $\mathrm{d}\log_{\mathbb{A}}$ for odd primes $p > 2$.

2.4. **Extending to animated prisms.** We note that the previous construction extends directly to the relative case: let $(A, I)$ be a transversal prism, considered as a base prism. Then the above construction carries out for prsims $(B^{\wedge}_{(p,I)}, IB^{\wedge}_{(p,I)})$ where $B$ is a free $\delta$-$A$-algebra $B$ generated by a finite set. By left deriving this construction, we get the (truncated) refined logarithm for every animated prisms over $(A, I)$. This is sufficient for the relative case when the base prism is transversal. However, there are two main disadvantages:

1. It is a priori unclear how much this construction does not depend on the base. More precisely, when we have two transversal prisms $(A, I)$ and $(B, J)$, and an animated prism $(C, K)$ over both $(A, I)$ and $(B, J)$, it is not obvious that the refined logarithm obtained by left deriving the previous construction does not depend on whether we choose $(A, I)$ or $(B, J)$ as the base prism.

2. Closely related, the construction does not extend directly to non-transversal prisms (and in particular, crystalline prisms). By [BL22, Prop 2.4.1], for every prism $(B, J)$ (and even for animated prisms, see Appendix A), there exists a transversal prism $(A, I)$ along with a map $(A, I) \to (B, J)$, but again as in the first disadvantage, it is unclear whether switching to $(A, I)$ would give rise to the "correct" answer.

Thus instead of left deriving the previous construction, we sketch how to "translate" the previous construction to the animated setting. The key is to find a derived counterpart of Lemma 2.1. Then the construction of Breuil–Kisin twists, norm maps, and the refined logarithm adapts.

First, it is technically convenient to introduce a derived version of *preprisms* in [BL22, Def 2.1.1].

**Definition 2.19.** *A* derived preprism *is a triple* $(A, I, A \leftarrow I)$ *(or a pair* $(A, I)$ *when there is no ambiguity), where*

- *$A$ is a derived $\delta$-ring (d'après [Hol23]). We will denote by $\varphi_A$ (or $\varphi$ if without ambiguity) the Frobenius lift $A \to A$.*
- *$(A, I)$ is a generalized Cartier divisor.*

*The $\infty$-category of derived preprisms is* $\mathrm{DAlg}^\delta \times_{\mathrm{DAlg}} \mathrm{GCart}$, *where* $\mathrm{GCart}$ *is the $\infty$-category of generalized Cartier divisors* $(A, I)$.

Recall that, for a $\delta$-ring $A$ and a non-zero-divisor $d \in A$, the generalized Cartier divisor $(A/d, (\delta(d)))$ does not depend on the orientation $d$ of the effective Cartier divisor $(d)$, cf. [Mao21, Lem 5.32]. In fact, this generalizes to all derived preprisms. To see this, we need the following technical lemma, which tells us when a map of modules factors through the multiplication by a generalized Cartier divisor.

**Lemma 2.20.** *Let $B$ be a derived ring, $(B, J)$ a generalized Cartier divisor, and $f : M \to N$ a map in $D(B)$. Then the datum of a factorization $M \xrightarrow{\exists} N \otimes_B^\mathbb{L} J \to N$ of $f$ in $D(B)$ is equivalent to the datum of a nullhomotopy of the map $M \otimes_B^\mathbb{L} (B/J) \to N \otimes_B^\mathbb{L} (B/J)$ in $D(B/J)$.*

*Proof.* Since there is a fiber sequence $N \otimes_B^\mathbb{L} J \to N \to N \otimes_B^\mathbb{L} (B/J)$, the datum of a factorization of $f$ as $M \xrightarrow{\exists} N \otimes_B^\mathbb{L} J \to N$ in $D(B)$ is equivalent to the datum of a nullhomotopy of the composite map $M \xrightarrow{f} N \to N \otimes_B^\mathbb{L} (B/J)$ in $D(B)$. Then the result follows from the adjunction $D(B) \xrightleftharpoons[]{\cdot \otimes_B^\mathbb{L} (B/\mathbb{L} J)} D(B/J)$. □

**Definition 2.21.** *Let $A$ be a derived ring, and $(A, I)$ a generalized Cartier divisor.*

- *A* null-homotopy *of the generalized Cartier divisor is simply a null-homotopy of the structure map $I \to A$ in $D(A)$.*
- *We say that $(A, I)$ is a* unit *if the structure map $I \to A$ is an equivalence in $D(A)$.*

**Notation 2.22.** *Let $A \to B$ be a map of derived rings, and $(A, I)$ a generalized Cartier divisor. Then we will denote by $(A, I) \otimes_A^\mathbb{L} B$ the base-changed generalized Cartier divisor $(B, I \otimes_A^\mathbb{L} B)$.*

**Notation 2.23.** *Let $A$ be a derived ring, $(A, I)$ a generalized Cartier divisor, and $(A, J)$ a derived Smith ideal. Then we will denote by $A/(I, J)$ the tensor product $(A/I) \otimes_A^\mathbb{L} (A/J)$ (which could be understood as a "Koszul complex" as well).*

We use this to factorize the generalized Cartier divisor $(A, \varphi^*(I)) \otimes_A^\mathbb{L} (A/I^p)$ through $(p)$.

*Construction* 2.24. Let $(A, I)$ be a derived preprism. We construct a generalized Cartier divisor $(A/I^p, \bar\delta_p(I))$ as follows. Note that we have canonical equivalences

$$\begin{aligned}
(A, \varphi^*(I)) \otimes_A^\mathbb{L} (A/^\mathbb{L} p) &= (A, I) \otimes_{A, \varphi}^\mathbb{L} A \otimes_A^\mathbb{L} (A/^\mathbb{L} p) \\
&= (A, I) \otimes_A^\mathbb{L} (A/^\mathbb{L} p) \otimes_{A/^\mathbb{L} p, \varphi}^\mathbb{L} (A/^\mathbb{L} p) \\
&= (A, I^p) \otimes_A^\mathbb{L} (A/^\mathbb{L} p)
\end{aligned}$$

of generalized Cartier divisors. This gives rise to a canonical null-homotopy of the generalized Cartier divisor $(A, \varphi^*(I)) \otimes_A^{\mathbb{L}} (A/(p, I^p))$. By Lemma 2.20, we get a map $(A, \varphi^*(I)) \otimes_A^{\mathbb{L}} (A/I^p) \to (A/I^p, p)$ of generalized Cartier divisors, i.e. a generalized Cartier divisor $(A/I^p, \bar{\delta}_p(I))$ such that the generalized Cartier divisor $(A, \varphi^*(I)) \otimes_A^{\mathbb{L}} (A/I^p)$ is the multiplication of generalized Cartier divisors $(A/I^p, p)$ and $(A/I^p, \bar{\delta}_p(I))$. We will denote by $(A/I, \bar{\delta}(I))$ the base-changed generalized Cartier divisor $(A/I^p, \bar{\delta}_p(I)) \otimes_{A/I^p}^{\mathbb{L}} (A/I)$.

**Definition 2.25.** *We say that a derived preprism $(A, I)$ is*
- *local if $(p, \pi_0(I))$ lies in the Jacobson radical $\mathrm{Rad}(\pi_0(A))$;*
- *distinguished if it is local and the generalized Cartier divisor $(A/I, \bar{\delta}(I))$ is unit;*
- *a derived prism if it is distinguished and $(p, I)$-complete.*

**Notation 2.26.** *Let $(A, I)$ be a derived preprism, and $r \in \mathbb{N}_{>0}$. We will denote by $I_r$ the generalized Cartier divisor $(A, I\varphi^*(I) \cdots (\varphi^{r-1})^*(I))$ being the multiplication of generalized Cartier divisors $(A, (\varphi^s)^*(I))$ for integers $0 \le s < r$, and by $J_r$ the generalized Cartier divisor $(A, I^p (\varphi^*(I) \cdots (\varphi^{r-1})^*(I))^{p-1} = I \cdot I_r^{p-1})$.*

We now generalize Lemma 2.1 to distinguished derived preprisms. For this, we first note that the multiplication of a Smith ideal by a generalized Cartier divisors can be analyzed by the following lemma.

**Lemma 2.27.** *Let $A$ be a derived commutative ring, $(A, I)$ a generalized Cartier divisor, and $(A, J)$ a derived Smith ideal. Then the canonical diagram*

$$\begin{array}{ccc} A/IJ & \longrightarrow & A/J \\ \downarrow & & \downarrow \\ A/I & \longrightarrow & A/(I, J) \end{array}$$

*of derived commutative $A$-algebras[2.2] is Cartesian.*

*Proof.* It suffices to check that this is a Cartesian diagram in $D(\mathbb{Z})$, which follows from comparing the fibers of horizontal maps. □

We now explain how to translate the proof of Lemma 2.1 to distinguished derived prisms. We start with the first part of the proof, namely, constructing $f_0 = p^{-1} f$ by invoking Lemma 2.20, just as what we did in Construction 2.24.

*Construction* 2.28. Let $(A, I)$ be a derived preprism. Then we construct a canonical map

$$(\varphi^r)^*(I) \otimes_A^{\mathbb{L}} (A/J_r) \longrightarrow p \, A/J_r$$

of generalized Cartier divisors of $A/J_r$ as follows. By Lemma 2.20, such a map is equivalent to the null-homotopy of the generalized Cartier divisor $(A, (\varphi^r)^*(I)) \otimes_A^{\mathbb{L}} (A/(p, J_r))$. Now we have canonical equivalences

$$\begin{aligned} (A, (\varphi^r)^*(I)) \otimes_A^{\mathbb{L}} (A/^{\mathbb{L}} p) &= (A, I^p) \otimes_A^{\mathbb{L}} (A/^{\mathbb{L}} p), \\ (J_r, A) \otimes_A^{\mathbb{L}} (A/^{\mathbb{L}} p) &= (A, I^{1+(p-1)(p^{r-1}+\cdots+1)}) \otimes_A^{\mathbb{L}} (A/^{\mathbb{L}} p) \\ &= (A, I^{p^r}) \otimes_A^{\mathbb{L}} (A/^{\mathbb{L}} p), \end{aligned}$$

---

[2.2]. We are mainly interested in the special case that $A$ is an animated ring, and $(A, J)$ is a connective derived Smith ideal. However, we stress that the diagram is Cartesian in the $\infty$-category of derived rings, which is stronger than being Cartesian in the $\infty$-category of animated rings.

which gives rise to the null-homotopy that we want.

The second part of the proof of Lemma 2.1 works verbatim. ∎

**Lemma 2.29.** *Let $(A, I)$ be a distinguished derived preprism. Then the map in Construction 2.28 is an equivalence.*

*Proof.* We can check this by passing to a pro-Zariski cover, thus without loss of generality, we may assume that the derived preprism $(A, I)$ is orientable, and we can pick an orientation $I = (d)$. Now it follows from the computation of $\varphi^r(d)$ in the proof of Lemma 2.1 (where we plug $d$ into $x$ there). □

The constructions of Breuil–Kisin twists, norm maps, and the refined logarithm then extend directly to animated prisms.

## 3. Prismatic Witt vectors and prismatic Hochschild homology

Kaledin introduced *polynomial Witt vectors* as a version of Witt vectors with coefficients. He used polynomial Witt vectors to define the *Hochschild–Witt homology*, as a non-commutative version of the de Rham–Witt complex. In this section, we introduce the prismatic analogues of both of them, discuss the relation to the prismatic Teichmüller map (Remark 2.11). Moreover, we indicate how to adapt the proof of Hesselholt's HKR theorem to establish an HKR-type theorem for prismatic Hochschild homology for $p$-completed polynomial algebras.

Let $(A, I)$ be a transversal prism. By [Sul23, Cons 3.3] (reviewed in Remark 2.4), we have a $\mathbb{T}_p$-Tambara functor $\underline{(A, I)}$ with $\underline{(A, I)}^{C_{p^{r-1}}} = A/I_r$ with trivial Weyl action.

**Definition 3.1.** *Let $(A, I)$ be a transversal prism, $r \in \mathbb{N}_{>0}$, and $M \in D(A/I)$ an $A/I$-module spectrum. The module $\mathbb{A}_r(M/A)$ of $r$-truncated prismatic Witt vectors is defined to be the $A/I_r$-module spectrum*
$$\left(M^{\otimes^{\mathbb{L}}_{\underline{(A,I)}} C_{p^{r-1}}}\right)^{C_{p^{r-1}}},$$
*where $(-)^{\otimes^{\mathbb{L}}_{\underline{(A,I)}} C_{p^{r-1}}} : D(A/I) \to \mathrm{Mod}_{\underline{(A,I)}}(\mathrm{Sp}^{gC_{p^{r-1}}})$ is the $C_{p^{r-1}}$-norm relative to the Tambara functor $\underline{(A, I)}$.*

We now relate prismatic Witt vectors to refined logarithm. First, it follows from a direct computation that

**Lemma 3.2.** *Let $G$ be a finite group, and $T$ a $G$-Tambara functor. Then the map*
$$T^e = \mathrm{End}_{D(T^e)}(T^e) \to \mathrm{End}_{D(T^G)}(T^G) = T^G$$
*induced by the composite functor*
$$D(T^e) \xrightarrow{N^G_e} \mathrm{Mod}_T(\mathrm{Sp}^G) \xrightarrow{(-)^G} D(T^G)$$
*coincides with the norm map $N : T^e \to T^G$.*

Applying this to the $C_{p^{r+1}}$-Tambara functor $\underline{(A, I)}$ (viewed as a $C_{p^s}$-Tambara functor for every $s \leq r + 1$), we get

**Corollary 3.3.** *Let $(A, I)$ be a transversal prism and $r \in \mathbb{N}_{>0}$. Then the map*
$$A/I = \mathrm{End}_{D(A/I)}(A/I) \longrightarrow \mathrm{End}_{D(A/I_r)}(A/I_r) = A/I_r$$

*induced by the r-truncated prismatic Witt vectors functor*

$$\mathbb{\Delta}_r \colon D(A/I) \longrightarrow D(A/I_r)$$

*coincides with the iterated norm map* $N_{r-1} \circ \cdots \circ N_1 \colon A/I \to A/I_r$ *in Construction 2.3.*

Another fact is that the $r$-truncated prismatic Witt vectors functor sends invertible modules to invertible modules, thus on invertible modules, it is the same as the map $H^1(\mathrm{Spf}(A/I); \mathbb{G}_m) \to H^1(\mathrm{Spf}(A/I_r); \mathbb{G}_m)$ induced by the iterated norm map $A/I \to A/I_r$. To see this, we first notice that the construction of the prismatic Witt vector functor is functorial in transversal prisms, by the rigidity of prisms. We just formulate a 1-categorical functoriality (which is sufficient to deduce the $\infty$-categorical functorial, as we can restrict to connective compact projective objects) as follows.

**Lemma 3.4.** *Let $(A, I) \to (B, K)$ be a map of transversal prisms. Then we have a commutative diagram*

$$\begin{array}{ccc} D(A/I) & \longrightarrow & D(A/I_r) \\ \downarrow & & \downarrow \\ D(B/K) & \longrightarrow & D(B/K_r) \end{array}$$

*of $\infty$-categories[3.1], where the horizontal arrows are $r$-truncated prismatic Witt vectors functors, and the vertical arrows are base changes.*

Note that, given an element $\bar{s} \in A/I$, the ($p$-completed) localization $A[s^{-1}]$ makes sense (by picking an arbitrary lift $s \in A$ of $\bar{s}$, since $I$ lies in the Jacobson radical $\mathrm{Rad}(A)$. Now we apply Lemma 3.4 to such localizations of transversal prisms, obtaining

**Lemma 3.5.** *Let $(A, I)$ be a transversal prism and $r \in \mathbb{N}_{>0}$. Then the $r$-truncated prismatic Witt vectors functor*

$$\mathbb{\Delta}_r \colon D(A/I) \longrightarrow D(A/I_r)$$

*sends an invertible $A/I$-module (in the heart) to an invertible $A/I_r$-module.*

Now we define the prismatic Hochschild homology. In fact, the prismatic Witt vectors give rise to a trace theory, and the prismatic Hochschild homology is the corresponding twisted Hochschild homology. In our case, the trace theory comes[3.2] from a Tambara functor, and in this case, the twisted Hochschild homology coincides with the relative Hochschild homology:

**Definition 3.6.** *Let $(A, I)$ be a transversal prism, and $\mathcal{C}$ a dualizable presentable stable $A/I$-linear $\infty$-category. The* prismatic Hochschild homology $\mathrm{HH}^{\mathbb{\Delta}}(\mathcal{C}/A) \in \mathrm{Mod}_{(A,I)}(\mathrm{Sp}^{g^{<\mathbb{T}}})$ *is defined to be the topological Hochschild homology of $\mathcal{C}$ relative to $(A, I)$ (cf. [ABG+18]), that is, by*

$$\mathrm{HH}^{\mathbb{\Delta}}(\mathcal{C}/A) := \mathrm{THH}(\mathcal{C}) \otimes^{\mathbb{L}}_{\mathrm{THH}(A/I)} (A, I),$$

*where the map $\mathrm{THH}(A/I) \to (A, I)$ of $\mathbb{T}$-normed rings is given by the universal property of $\mathrm{THH}$ (cf. [ABG+18]).*

---

3.1. Although the categories in the commutative diagram are stable, the horizontal functors are not exact, thus it is not a commutative diagram of stable $\infty$-categories.

3.2. In our forthcoming work [Mao], we will explain how to produce a trace theory from a normed category.

*Example* 3.7. Let $S$ be a $p$-torsion-free perfectoid ring, and $\mathcal{C}$ a dualizable presentable stable $S$-linear $\infty$-category. Let $(A, I)$ denote the prism $(A_{\inf}(S), (\xi))$ corresponding to the perfectoid ring $S$. Then the prismatic Hochschild homology $\mathrm{HH}^{\mathbb{\Delta}}(\mathcal{C}/A)$ coincides with the relative Hochschild–Witt homology $\mathrm{HH}(\mathcal{C}/S)$ in our forthcoming work [Mao].

We now formulate a conjectural HKR type theorem for the prismatic Hochschild homology. For this, we first have to introduce the *prismatic de Rham complex*. Such a complex was implicitly considered in [Mol20, Prop 4.10] when the base prism is oriented (also compare with $\mathcal{W}_r^n(D)$ in [BMS18, §11.1.2]).

*Construction* 3.8. Let $(A, I)$ be a transversal prism, $r \in \mathbb{N}_{>0}$, and $R$ a $p$-completely smooth $A/I$-algebra. We construct the $r$-truncated *prismatic de Rham complex* $(\mathbb{\Delta}_r \Omega_{R/A}^*)$ as follows. For every $n \in \mathbb{Z}$, we define

$$\mathbb{\Delta}_r \Omega_{R/A}^n := H^n\Big(\mathbb{\Delta}_{R/A} \otimes_A^{\mathbb{L}} (I_r/I_r^2)^{\otimes_{A/I_r}^{\mathbb{L}} n}\Big),$$

with differentials

$$\mathrm{d} : \mathbb{\Delta}_r \Omega_{R/A}^n \longrightarrow \mathbb{\Delta}_r \Omega_{R/A}^{n+1}$$

given by applying $H^n((-)\{n\})$ to the Bockstein map

$$\mathbb{\Delta}_{R/A} \otimes_A^{\mathbb{L}} (A/I_r) \longrightarrow \mathbb{\Delta}_{R/A} \otimes_A^{\mathbb{L}} (I_r/I_r^2)[1].$$

**Conjecture 3.9.** *Let $(A, I)$ be a transversal prism, $r \in \mathbb{N}_{>0}$, and $R$ a $p$-completely smooth $A/I$-algebra. Recall that the $\mathbb{T}$-action on $\mathrm{HH}^{\mathbb{\Delta}}(R/A)^{C_{p^{r-1}}}$ gives rise to a cochain complex $\big(\pi_* \mathrm{HH}^{\mathbb{\Delta}}(R/A)^{C_{p^{r-1}}}, \mathrm{d}_{\mathrm{HH}}\big)$. Then we have an isomorphism*

$$\big(\pi_* \mathrm{HH}^{\mathbb{\Delta}}(R/A)^{C_{p^{r-1}}}, \mathrm{d}_{\mathrm{HH}}\big) \cong (\mathbb{\Delta}_r \Omega_{R/A}^*, \mathrm{d})$$

*of cochain complexes, which is functorial with respect to the data $(A, I, R)$.*

We could produce a non-functorial comparison isomorphism in Conjecture 3.9 for $p$-completed polynomial $A/I$-algebras $R$. In fact, it is functorial with respect to permutations of chosen generators. In fact, this is a consequence of the proof of Hesselholt's HKR theorem, which we now sketch. We are unable to give any insight on this argument, thus we omit the combinatorial details.

**Notation 3.10.** *Let $R$ be an $\mathbb{E}_\infty$-ring, and $S$ a finite set. Then we will denote by $R[x_s \,|\, s \in S]$ the flat polynomial $R$-algebra generated by $S$.*

*Remark* 3.11. Let $S$ be a finite set. Then we have

$$\mathrm{THH}(\mathbb{S}[x_s \,|\, s \in S]) = \Sigma_{\mathbb{T}}^\infty [B^{\mathrm{cy}}(\mathbb{N}[S])]_+$$

as cyclonic spectra, where $\mathbb{N}[S]$ is the free monoid generated by $S$. Let $T$ be a $\mathbb{T}$-Tambara functor. Then

$$\begin{aligned} \mathrm{THH}(T^e[x_s \,|\, s \in \mathbb{S}]) \otimes_{\mathrm{THH}(T^e)} T &= (\mathrm{THH}(T^e) \otimes \mathrm{THH}(\mathbb{S}[x_s \,|\, s \in S])) \otimes_{\mathrm{THH}(T^e)} T \\ &= T \otimes \Sigma_{\mathbb{T}}^\infty [B^{\mathrm{cy}}(\mathbb{N}[S])]_+ \end{aligned}$$

where $B^{\mathrm{cy}}(-)$ is the cyclic bar construction on $\mathbb{E}_1$-monoids.

We first review the description of the $\mathbb{T}$-anima $B^{\mathrm{cy}}(\mathbb{N}[S])$ in [Hes96, §2.2]. Let $\mathbb{T}(0)$ denote the point with trivial $\mathbb{T}$-action. For every $l \in \mathbb{Z} \setminus 0$, let $\mathbb{T}(l)$ denote the $\mathbb{T}$-anima $\rho_l^* \mathbb{T}$, where $\rho_l : \mathbb{T} \to \mathbb{T}$ corresponds to $l \in \mathbb{Z} = \pi_1(\mathbb{T})$ (represented by $(-)^l : S^1 \to S^1$). Then the $\mathbb{T}$-anima $B^{\mathrm{cy}} \mathbb{N}$ can be identfied with the disjoint union $\coprod_{l \in \mathbb{N}} \mathbb{T}(l)$.

Now for every map $l: S \to \mathbb{Z}$, let $\mathbb{T}(l)$ denote the product $\prod_{s \in S} \mathbb{T}(l(s))$. Then the $\mathbb{T}$-anima $B^{\mathrm{cy}}(\mathbb{N}[S])$ is equivalent to the disjoint union $\coprod_{l: S \to \mathbb{N}} \mathbb{T}(l)$. We will also denote $v(l) := \min \{v_p(l(S))\}$.

The proof of [Hes96, Prop 2.2.5] implies that, for every $C_{p^r}$-spectrum $T$, and every map $l: S \to \mathbb{N}$, up to $p$-localization, the spectrum

$$(T \otimes \Sigma^\infty_\mathbb{T} \mathbb{T}(l)_+)^{C_{p^r}}$$

can be identfied with $T^{C_{p^w}} \otimes (\mathbb{T} \times \cdots \times \mathbb{T})$, where $w := \min \{r, v(l)\}$. The argument in [Hes96, §2.3] essentially implies that

**Proposition 3.12.** *Let $T$ be a $\mathbb{T}$-Mackey functor, $r \in \mathbb{N}_{>0}$, and $l: S \to \mathbb{N}$ a map. Let $a: S \to \mathbb{N}$ be the map $s \mapsto l(s) - r$. Then up to $p$-completion, for every $n \in \mathbb{N}$, the $n$-th homotopy group*

$$\pi_n\bigl((T \otimes \Sigma^\infty_\mathbb{T} \mathbb{T}(l)_+)^{C_{p^{r-1}}}\bigr)$$

*can be identified with*

$$\bigoplus_{(I_0, \ldots, I_n) \in P_a} T^{C_{p^{r-1-u(a)}}},$$

*where $u(a) := \max \{-v(a), 0\}$, and the set $P_a$ is as defined in [BMS18, §10.4].*

Now let $(A, I)$ be a transveral prism, and $S$ a finite set. Then [Mol20, Prop 4.9] shows that the $n$-th relative prismatic cohomology $H^n(\mathbb{\Delta}_{(A/I)\langle x_s \mid s \in S \rangle / A})$ of the $p$-completed polynomial $A/I$-algebra $(A/I)\langle x_s \mid s \in S \rangle$ has the same form.

## Appendix A. Transversal approximation of derived prisms

In this appendix, we briefly explain that the transversal approximation in [BL22, §2.4] works well for derived prisms. First, every derived prism admits a connective cover, thus we may replace "derived" by "animated".

The following smooth spreadout works verbatim.

**Lemma A.1.** ([BL22, Lem 2.4.2]) *Let $R$ be an animated ring, and $M$ a finite projective $R$-module. Then there exist a smooth $\mathbb{Z}$-algebra $R_0$, a finite projective $R_0$-module $M_0$, a map $R_0 \to R$ of animated rings, and an isomorphism $M \simeq M_0 \otimes^\mathbb{L}_{R_0} R$.*

To adapt the proof of [BL22, Prop 2.4.1], one has to establish a slightly stronger universal property than [BL22, Prop 2.1.10].

**Lemma A.2.** *The fully faithful functor*

$$\{\text{animated prisms}\} \hookrightarrow \{\text{animated preprisms}\}$$

*admits a left adjoint. Let $(A, I)$ be an animated preprism, and $(B, J)$ the image of $(A, I)$ under the left adjoint above. Then the map $A \to B$ underlying the unit natural transformation $(A, I) \to (B, J)$ is $(p, I)$-completely flat. In particular, if $(A, I)$ is a transversal preprism, then $(B, J)$ is a transversal prism.*

*Proof.* Let $(A, I)$ be an animated preprism. We work internally in the $(p, I)$-complete derived category $D(A)^\wedge_{(p,I)}$. Let $\mathcal{C}$ be the site $\{A[f^{-1}] \mid f \in \pi_0(A)\}^{\mathrm{op}}$ of basic Zariski open subsets of $\mathrm{Spf}(\pi_0(A))$. Then the presheaves

$$C \longmapsto \{(p, I)\text{-complete animated preprisms over } (C, IC)\}, \text{ and}$$
$$C \longmapsto \{\text{animated prisms over } (C, IC)\}$$

of $\infty$-categories are sheaves, thus it suffices to reduce to the case that $I = (d)$ is principal, in which case $B = A[\delta(d)^{-1}]$ is $A$-flat. Then the result follows from the Zariski descent. $\square$

Then the proof in [BL22, §2.4] leads to

**Proposition A.3. ([BL22, Prop 2.4.1])** *Let $(A, I)$ be an animated prism. Then there exists a transversal prism $(B, J)$ along with a map $(B, J) \to (A, I)$ of animated prisms, such that $|B| \leq 2^{\aleph_0}$.*

**Proposition A.4. ([BL22, Prop 2.4.5 & Rem 2.4.6])** *Let $(A, I)$ be a transversal prism. Then for every animated prism $(B, J)$, the canonical map $(B, J) \to (A, I) \amalg (B, J)$ is flat.*

*Proof.* Note that the pushout $(A, I) \amalg (B, J)$ is the animated prismatic envelope of $A \otimes_{\mathbb{Z}}^{\mathbb{L}} B \to (A/I) \otimes_{\mathbb{Z}}^{\mathbb{L}} (B/J)$ over the animated prism $(B, J)$. The result follows immediately from the Hodge–Tate comparison of the animated prismatic envelope, cf. [Mao21, §5.3]. $\square$

**Corollary A.5. ([BL22, Cor 2.4.7])** *Let $\kappa$ be a regular cardinal. Then any $\kappa$-small coproduct of transversal prisms, taken in the $\infty$-category of animated prisms, is also transversal.*

of $\infty$-categories are sheaves, thus it suffices to reduce to the case that $I = (d)$ is principal, in which case $B = A[\delta(d)^{-1}]$ is $A$-flat. Then the result follows from the Zariski descent. $\square$

Then the proof in [BL22, §2.4] leads to

**Proposition A.3. ([BL22, Prop 2.4.1])** *Let $(A, I)$ be an animated prism. Then there exists a transversal prism $(B, J)$ along with a map $(B, J) \to (A, I)$ of animated prisms, such that $|B| \leq 2^{\aleph_0}$.*

**Proposition A.4. ([BL22, Prop 2.4.5 & Rem 2.4.6])** *Let $(A, I)$ be a transversal prism. Then for every animated prism $(B, J)$, the canonical map $(B, J) \to (A, I) \amalg (B, J)$ is flat.*

*Proof.* Note that the pushout $(A, I) \amalg (B, J)$ is the animated prismatic envelope of $A \otimes_{\mathbb{Z}}^{\mathbb{L}} B \to (A/I) \otimes_{\mathbb{Z}}^{\mathbb{L}} (B/J)$ over the animated prism $(B, J)$. The result follows immediately from the Hodge–Tate comparison of the animated prismatic envelope, cf. [Mao21, §5.3]. $\square$

**Corollary A.5. ([BL22, Cor 2.4.7])** *Let $\kappa$ be a regular cardinal. Then any $\kappa$-small coproduct of transversal prisms, taken in the $\infty$-category of animated prisms, is also transversal.*

## Bibliography

of $\infty$-categories are sheaves, thus it suffices to reduce to the case that $I = (d)$ is principal, in which case $B = A[\delta(d)^{-1}]$ is $A$-flat. Then the result follows from the Zariski descent. $\square$

Then the proof in [BL22, §2.4] leads to

**Proposition A.3. ([BL22, Prop 2.4.1])** *Let $(A, I)$ be an animated prism. Then there exists a transversal prism $(B, J)$ along with a map $(B, J) \to (A, I)$ of animated prisms, such that $|B| \leq 2^{\aleph_0}$.*

**Proposition A.4. ([BL22, Prop 2.4.5 & Rem 2.4.6])** *Let $(A, I)$ be a transversal prism. Then for every animated prism $(B, J)$, the canonical map $(B, J) \to (A, I) \amalg (B, J)$ is flat.*

*Proof.* Note that the pushout $(A, I) \amalg (B, J)$ is the animated prismatic envelope of $A \otimes_{\mathbb{Z}}^{\mathbb{L}} B \to (A/I) \otimes_{\mathbb{Z}}^{\mathbb{L}} (B/J)$ over the animated prism $(B, J)$. The result follows immediately from the Hodge–Tate comparison of the animated prismatic envelope, cf. [Mao21, §5.3]. $\square$

**Corollary A.5. ([BL22, Cor 2.4.7])** *Let $\kappa$ be a regular cardinal. Then any $\kappa$-small coproduct of transversal prisms, taken in the $\infty$-category of animated prisms, is also transversal.*

## Bibliography

[ABG+18] Vigleik Angeltveit, Andrew J. Blumberg, Teena Gerhardt, Michael A. Hill, Tyler Lawson, and Michael A. Mandell. Topological cyclic homology via the norm. *Doc. Math.*, 23:2101–2163, 2018.

[Bha23] Bhargav Bhatt. Crystals and Chern classes. *ArXiv e-prints*, oct 2023.

[BL22] Bhargav Bhatt and Jacob Lurie. Absolute prismatic cohomology. *ArXiv e-prints*, jan 2022.

[BMS18] Bhargav Bhatt, Matthew Morrow, and Peter Scholze. Integral $p$-adic Hodge theory. *Publications Mathématiques. Institut de Hautes Études Scientifiques*, 128:219–397, 2018.

[BMS19] Bhargav Bhatt, Matthew Morrow, and Peter Scholze. Topological Hochschild homology and integral $p$-adic Hodge theory. *Publications Mathématiques. Institut de Hautes Études Scientifiques*, 129:199–310, 2019.

[DKNP23] Emanuele Dotto, Achim Krause, Thomas Nikolaus, and Irakli Patchkoria. Witt vectors with coefficients and TR. *ArXiv e-prints*, dec 2023.

[Hes96] Lars Hesselholt. On the $p$-typical curves in Quillen's $K$-theory. *Acta Math.*, 177(1):1–53, 1996.

[Hol23] Adam Holeman. Derived $\delta$-Rings and Relative Prismatic Cohomology. *ArXiv e-prints*, mar 2023.

[Lur18] Jacob Lurie. Spectral algebraic geometry. https://www.math.ias.edu/~lurie/papers/SAG-rootfile.pdf, feb 2018.

[Mao] Zhouhang Mao. Polynomial Witt vectors and Hochschild–Witt homology via norm. In preparation.

[Mao21] Zhouhang Mao. Revisiting derived crystalline cohomology. *ArXiv e-prints*, jul 2021.

[Mol20] Semen Molokov. Prismatic cohomology and de Rham-Witt forms. *ArXiv e-prints*, aug 2020.

[Rea23] Thomas Read. $G$-typical Witt vectors with coefficients and the norm. *ArXiv e-prints*, may 2023.

[Sul23] Yuri J. F. Sulyma. Prisms and Tambara functors I: Twisted powers, transversality, and the perfect sandwich. *ArXiv e-prints*, sep 2023.